\documentstyle[12pt]{article}
\begin{document}
\parindent=45pt
\def\z{{\bf Z}}
\def\n{{\bf N}}
\def\c{{\bf C}}
\thispagestyle{empty}
\ \vskip 30mm
{\bf ON THE POLAR DECOMPOSITION OF CIRCULAR VARIABLES}
\vskip 30pt
Teodor Banica
\vskip 30pt
We find an elementary proof for Voiculescu's theorem on the polar decomposition of circular
variables.
\vskip 30pt

\baselineskip=18pt

{\bf 0. INTRODUCTION: }In non commutative probability theory (see [1]) an important role is played
by the Haar-unitaries, which appear naturally in the von Neumann algebras of free groups
$W^*({\z}^{*n})$, and by the circular variables, which appear naturally in the algebras of creation
operators on the full Fock spaces $W^*({\n}^{*n})$.

In [2] Voiculescu finds the polar decomposition of the circular variables. In order to work at the
same time with circular variables and with Haar-unitaries, he uses approximation by random
Gaussian matrices.

Following an idea of G.Skandalis, we give in this note an elementary proof of this theorem: it
turns out that the result is an immediate consequence of some combinatorial properties of the
monoid $\z *\n $. 

In what follows, we introduce (see 4) a certain class of monoids, which contains $\n$, all the
groups, and is stable under free products (see 8). Consequently, for such monoids $M$, many
variables having interesting $*$-distributions (semicircular and circular variables, Haar-unitaries)
appear naturally in the algebras $W^*(M)$.

For such monoids, some combinatorial properties (see 5 and 7) allow us to find (and easily
manipulate) ``many'' circular systems in their algebras $W^*(M)$ (see 6.2). The polar
decomposition of the circular variables follows easily.

{\bf 1. NOTATION: }In what follows, we will denote (by abuse of language) by {\bf ``monoid''} a
{\bf countable unital monoid, which can be embedded in a group}. For such a monoid $M$, the symbol
$M^{-1}$ will denote, when there are no confusions, the subset $\{ m^{-1},m\in M\}$ of some group
containing $M$.

{\bf 2. DEFINITIONS: }Let $M$ be a monoid and $l^2(M)$ the Hilbert space of square summable
functions from $M$ to ${\bf C}$, with $({\delta}_m)_{m\in M}$ the canonical orthonormal basis.
Using the left simplifiability of $M$ one can define, as for discrete groups, the embedding of
monoids $(M,\cdot )\rightarrow (B(l^2(M)),\circ )$ by $\lambda _{M} (m)\delta _n=\delta
_{mn}$. Let $W^*(M)$ be the Von Neumann algebra generated by ${\lambda}_M(M)$. Together with the
canonical state ${\tau}_M(T)=<T{\delta}_e,{\delta}_e>$ it is a non commutative $W^*$-probability
space.

{\bf 3.1. REMARK: }The operators in $\lambda _M(M)$ are isometries, but not necessarely unitaries,
as in the group case. Indeed, for every $m\in M$, ${\lambda}_M(m)^*$ is given by
${\lambda}_M(m)^*({\delta}_n)={\Sigma}_{x\in
M}<{\lambda}_M(m)^*{\delta}_n,{\delta}_x>{\delta}_x={\Sigma}_{x\in M}{\delta}_{n ,mx}{\delta}_x$,
so that ${\lambda}_M(m)^*{\lambda}_M(m)=1$.

{\bf 3.2. REMARK: }It is easy to see that $l^2(\n^{*I})$ is the full Fock space over ${\c}^I$. By
this identification, $(W^*({\n}^{*I}),{\tau}_{{\n}^{*I}})$ is the algebra of creation operators,
with the canonical state associated to the vacuum vector.

{\bf 3.3. REMARK: }Let $M\subset N$ be monoids (so $l^2(M)\subset l^2(N)$). For $m,\,m^{\prime}\in
M$ one has ${\lambda}_M(m){\delta}_{m^{\prime}}={\lambda}_N(m){\delta}_{m^{\prime}}$, so if we
suppose $M(N-M)=N-M$ then ${\lambda}_M(m)^*{\delta}_{m^{\prime}}={\Sigma}_{x\in
M}{\delta}_{m^{\prime},mx}{\delta}_x={\Sigma}_{x\in
N}{\delta}_{m^{\prime},mx}{\delta}_x={\lambda}_N(m)^*{\delta}_{m^{\prime}}$. In particular, if
$m_1...m_k\in M$ and ${\alpha}_1...{\alpha}_k$ are exponents $\in \{1,*\}$ then
${\lambda}_M(m_1)^{{\alpha}_1}...{\lambda}_M(m_k)^{{\alpha}_k}{\delta}_e={\lambda}_N(m_1)^
{{\alpha}_1}...{\lambda}_N(m_k)^{{\alpha}_k}{\delta}_e$.

It follows that if $M\subset N$ are monoids such that $M(N-M)=N-M$ then for every family $\{
a_i\}_{i\in I}$ of elements in $M$, the $*$-distribution joint to $\{\lambda_N(a_i)\}_{i\in I}$ is
equal to the $*$-distribution joint to $\{ \lambda_M(a_i)\}_{i\in I}$.

{\bf 3.4. REMARK: }Let $\prod _{i\in I}M_i$ be a direct product of monoids and $a\in M_k\,,\,k\in
I$. Then, by (3.3) the $*$-distribution of $\lambda _{M_k}(a)\in W^*(M_k )$ is equal to the
$*$-distribution of $\lambda _{\prod M_i}(a)\in W^*(\prod M_i)$. Moreover, let $b\in M_j$ with
$j\neq k$. Then it is easy to see that $\lambda _{\prod M_i}(a)$ and $\lambda _{\prod M_i}(b)$ are
independent.

{\bf 3.5. REMARK: }Let $\ast _{i\in I}M_i$ be a free product of monoids and $a\in M_k\,,\,k\in I$.
Then, by (3.3) the $*$-distribution of $\lambda _{M_k}(a)\in W^*(M_k )$ is equal to the
$*$-distribution of $\lambda _{*M_i}(a)\in W^*(*M_i)$.

\vskip 2mm We are now interested to find semicircular variables in the algebras of monoids. Let us
introduce some definitions related to the combinatorics of free monoids:

{\bf 4. DEFINITION: }Let $N$ be a monoid. Consider the following order relation on it:
$a{\preceq}_Nb$ if and only if $b\in aN$. We say that $N$ is in the class $E$ if it satisfies one
of the (obvious) equivalent conditions:

(4.1.) for ${\preceq}_N$ every bounded subset is totally ordered.

(4.2.) $(a \preceq c , b \preceq c \Rightarrow a \preceq b$ or $b \preceq a$).

(4.3.) $aN\cap bN\neq \emptyset \Rightarrow aN\subset bN$ or $bN\subset aN$.

(4.4.) $NN^{-1}\cap N^{-1}N=N\cup N^{-1}$.

{\bf 5. DEFINITION: }Let $(a_i)_{i\in I}$ be a family of elements in a monoid $N$.\\
We call it ``code'' if it satisfies the following conditions:

(5.1.) the monoid $M$ generated by the $a_i$'s is isomorphic to  ${\n}^{*I}$ by
$a_i\mapsto e_i$.

(5.2.) $M(N-M)=N-M$.\\
We call it ``prefix'' if it satisfies the following condition:

(5.3.) $a_i\in a_jN\Rightarrow i=j$ (ie. the $a_i$'s are not comparable by ${\preceq}_N$).

{\bf 6.1. REMARK: }These are extensions of the classical notions of code and prefix, which already
appeared in the combinatorial theory of free monoids. A well known result (see [3]) asserts that on
free monoids the prefixes are the codes. We will extend this result to all monoids in the class $E$.

{\bf 6.2. REMARK: }Let $(a_i,b_i)_{i\in I}$ be a code. Then by (3.2) and (3.3), the family
$({\lambda}_N(a_i),{\lambda}_N(b_i))_{i\in I}$ has the same $*$-distribution as a family of
creation operators associated to a family of $2I$ orthonormal vectors, acting on the Fock space. In
particular $(1/2({\lambda}_N(a_i)+{\lambda}_N(b_i)^*))_{i\in I}$ is a circular family. Thus the
following proposition is a nice criterion for finding circular systems in the algebras $W^*(N)$ of
monoids (in the class $E$).

{\bf 7. PROPOSITION:} {\it For a monoid $N\in E$, a family $(a_i)_{i\in I}\subset N$ having at
least two elements is a prefix if and only if it is a code.}

{\bf PROOF: }Let $(a_i)_{i\in I}$ be a code which is not a prefix. Suppose for instance that
$a_i=a_jn$ with $i\neq j, n\in N$. By (5.2) $n$ is in the monoid $M$ generated by the $a_k$'s and
$a_i=a_jn$ with $i\neq j$, so $M$ cannot be free, contradiction.

Suppose now that $(a_i)_{i\in I}$ is a prefix and let
$A=a_{i_1}^{{\alpha}_1}...a_{i_n}^{{\alpha}_n}m=a_{j_ 1}^{{\beta}_1}...a_{j_s}^{ {\beta}_s}$ with
$m\in N$. One has $a_{i_1}\preceq A$, $a _{j_1}\preceq A$, so by (4.2) and (5.3), $i_1 =j_1$. We
simplify $A$ to the left by $a_{i_1}$ (recall that all the monoids we consider are bisimplifiable);
a reccurence on $\Sigma {\alpha}_ i$ shows that $n\leq s,\,a_{i_k}=a_{j_k}\,(\forall k\leq
n),\,{\alpha}_k={\beta}_k\,(\forall k<n),\,{\alpha}_n\leq
{\beta}_n,\,m=a_{j_n}^{{\beta}_n-{\alpha}_n}a_{j_{n+1}}^{{\beta}_{n+1}}...a_{j_s
}^{{\beta}_s}$. Finally, $m$ is in the monoid generated by the $a_i$'s, so (5.2) is true. Moreover,
for $m=e$ we obtain $n=s,\,a_{j_k}=a_{i_k},\,{\alpha}_k={\beta}_k,\,(\forall k\leq n )$ so
the $a_i$'s generate freely $M$ and $(a_i)_{i\in I}$ is a code.

{\bf 8. PROPOSITION:} {\it (8.1.) all the groups are in $E$.\\
(8.2.) the positive parts of totally ordered abelian groups are in $E$.\\
(8.3.) if $G$ is a group and $M\in E$, then $M\times G\in E$.\\
(8.4.) if $A_1$, $A_2$ are in $E$, then the free product $A_1*A_2$ is in $E$.}

{\bf PROOF: }(8.1) et (8.2) are obvious ($M$ is totally ordered by ${\preceq}_M$).

Remark: Reciprocally, if $M$ is an abelian monoid in  $E$, then one can easily construct a total
order on its Grothendieck group $K(M)$ such that $M=\{ g\in K(M), g\succeq 0\}$.

Let $G$ be a group and $M\in E$. Using (4.4) we have $(M\times G)(M\times G)^{-1}\cap (M\times
G)^{-1}(M\times G)=(M\times G)(M^{-1}\times G)\cap (M^{-1}\times G)(M\times G)=(MM^{-1}\times
G)\cap (M^{-1}M\times G)=(MM^{-1}\cap M^{-1}M)\times G=(M\cup M^{-1})\times G=(M
\times G)\cup (M^{-1}\times G)=(M\times G)\cup (M\times G)^{-1 }$, so we proved (8.3).

We prove now (8.4). Let $a,b,c\in A_1*A_2$ such that $ab=c$.

Write $a=x_1...x_n$, $b=y_1...y_m$, $c=z_1...z_p$ as reduced words. Let $s$ be such that
$x_ny_1=1$,...,$x_{n-s+1}y_s=1$, but $x_{n-s}y_{s+1}\neq 1$. Let
$u=x_{n-s+1}...x_n=(y_1...y_s)^{-1}$. Then $c=ab=x_1...x_{n-s}y_{s+1}...y_m$. Let $i\in\{ 1,2\}$ be
such that $z_{n-s}\in A_i$. There are two cases:

- if $x_{n-s}\in A_1$ and $y_{s+1}\in A_2$ or if $x_{n-s}\in A_2$ and $y_{s+1}\in A_1$, then
$x_1...x_{n-s}y_{s+1}...y_m$ is a reduced word. In particular, $x_1=z_1$,...,$x_{n-s}=z_{n-s}$.
Thus $a=z_1...z_{n-s}u$ with $u$ invertible.

- if $x_{n-s},y_{s+1}\in A_i$ then $x_1=z_1$,...,$x_{n-s-1}=z_{n-s-1}$ and $x_{n-s}y_{s+1}=z_{n-s}$.
In this case $a=z_1...z_{n-s-1}x_{n-s}u$ with $u$ invertible.

Remark that in both cases we obtained that $a$ is of the form $z_1...z_fxu$ for some $f$, with $u$
invertible and such that if $z_{f+1}\in A_i$, then there exists $y\in A_i$ with $xy=z_{f+1}$ (take
$f=n-s-1$ and $x=z_{n-s},y=1$ in the first case, $x=x_{n-s},y=y_{s+1}$ in the second one).

Suppose now that $A_1,A_2\in E$ and let $a,b,a^{\prime},b^{\prime}\in A_1*A_2$ such that
$ab=a^{\prime}b^{\prime}$. Let $z_1...z_p$ be the decomposition of $ab=a^{\prime}b^{\prime}$ as a
reduced word. Then we can decompose $a=z_1...z_fxu$ and
$a^{\prime}=z_1...z_{f^{\prime}}x^{\prime}u^{\prime}$ as above. We have to show that
$a=a^{\prime}m$ or that $a^{\prime}=am$ for some $m\in A_1*A_2$. There are three cases:

- if $f<f^{\prime}$, then $a^{\prime}=au^{-1}yz_{f+2}...z_{f^{\prime}}x^{\prime}u^{\prime}$.

- if $f^{\prime}<f$, then $a=a^{\prime}u^{\prime -1}z_{f^{\prime}+2}...z_fxu$.

- if $f=f^{\prime}$, then $xy=x^{\prime}y^{\prime}=z_{f+1}\in A_i$ for some $i\in\{ 1,2\}$. As
$A_i\in E$, we have that $x=x^{\prime}m$ or $x^{\prime}=xm$ for some $m\in
A_i$, so that $a=a^{\prime}u^{\prime -1}mu$ or $a^{\prime}=aumu^{\prime}$.

The proof of (8.4) is now complete.

{\bf 9. PROPOSITION: }{\it (9.1.) Let $M\subset N$ be two monoids in the class $E$ such that
$M(N-M)=N-M$. Let $\lambda=\lambda_N$. Then every element $x$ of the $*$-algebra generated by 
$\lambda (M)$ could be written as $x=\Sigma a_i\lambda(p_i)\lambda (q_i)^*$, with $p_i,q_i\in M$.

(9.2.) Let $A,B\in E$, $\lambda=\lambda_{A*B}$, $\tau=\tau_{A*B}$, and let $x$ be an element of the
$*$-algebra generated by $\lambda (A)$ such that $\tau (x)=0$. Denote by $W_A$ the set of reduced
words beginning by an element of $A$, and by $W_B$ the set of reduced words beginning by an
element of $B$. Then $x$ maps $l^2(W_B\cup \{ e\} )$ into $l^2(W_A)$.

(9.3.) Let $A,B\in E$. Then $\lambda _{A*B}(A)$ and $\lambda _{A*B}(B)$ are $*$-free.}

{\bf PROOF: }It is enough to prove (9.1) for $x=\lambda( m)^*\lambda (n)$ with $m,n\in M$; the
general case will follow easily. Remark that $x=\lambda( m)^*\lambda (n)$ is different from $0$ iff
$\exists a,b\in N$ such that $<\lambda( m)^*\lambda (n)\delta_a,\delta_b>\neq 0$, ie. if $na=mb$.
By (4.2), $\exists c\in N$ with $n=mc$ or with $m=nc$. Moreover, as $M(N-M)=N-M$, it follows that
$c\in M$. Thus $x=\lambda( m)^*\lambda (n)\neq 0\Rightarrow x=\lambda (c)$ or $x=\lambda (c)^*$
with $c\in M$, and this finishes the proof.

For proving (9.2), we apply (9.1) with $M=A$ and $N=A*B$ for writing $x=\Sigma
a_i\lambda(p_i)\lambda (q_i)^*$, with $p_i,q_i\in A$. Remark that $\tau (\lambda(p_i)\lambda
(q_i)^*)=\Sigma_x\delta_{e,p_ix}\delta_{e,q_ix}$ is nonzero iff $p_i=q_i$=invertible, and in this
case $\lambda(p_i)\lambda (q_i)^*=1$. As $\tau (x)=0$, it follows that we may write $x=\Sigma
a_i\lambda(p_i)\lambda (q_i)^*$, such that $\tau (\lambda(p_i)\lambda (q_i)^*)=0$ for every $i$. By
linearity, it is enough to prove (9.2) for $x=\lambda(p_i)\lambda (q_i)^*$.

Let $m\in  W_B\cup\{ e\}$ and suppose that $x\delta_m\neq 0$. Then $\lambda (q_i)^*\delta_m\neq
0$ implies that $m=q_ic$ for some word $c\in A*B$. As $q_i\in A$ and $m\in  W_B\cup\{ e\}$, it
follows that $q_i$ is invertible. In this case, $x\delta_m=\delta_{p_iq_i^{-1}m}\in l^2(W_A)$
(recall that $p_iq_i^{-1}=1\Rightarrow\tau (x)=1$).

Finally, (9.3) follows from (9.2). Indeed, let $P=x_n...x_1$ be a product of elements in $ker(\tau
)$, such that $x_{2k}$ is in the $*$-algebra generated by $\lambda (B)$ and $x_{2k+1}$ is in the
$*$-algebra generated by $\lambda (A)$. Then $x_1\delta_e\in l^2(W_A)$, so that $x_2x_1\delta_e\in
l^2(W_B)$ etc. By a reccurence, $P\delta_e$ is in $l^2(W_A)$ or in $l^2(W_B)$, and this implies
that $\tau (P)=0$.

{\bf 10. PROPOSITION: }{\it Consider (in some non commutative probability space) a Haar-unitary $u$,
$*$-free from a semicircular $s$. Then $us$ is a circular variable.}

{\bf PROOF: }
Denote by $z$ the image of $1\in \z$ and by $n$ the image of $1\in \n$ by the canonical embeddings
into the free product $\z *\n$. Let $\lambda ={\lambda} _{\z *\n}$. By (8), $\z *\n\in E$.
$(zn,nz^{-1})$ is obviously a prefix, so by the criterion (7), it is a code. By (6.2), $1/2(\lambda
(zn)+\lambda (nz^{-1})^*)$ is circular. But $1/2(\lambda (zn)+\lambda (nz^{-1})^*)=us$ where:

- $u:={\lambda}(z)$ is a Haar-unitary (see (3.5)).

- $s:=1/2(\lambda (n)+\lambda (n)^*)$ is semicircular (see (3.5) and (3.2)).

- $u$ and $s$ are $*$-free (by (9.3)). 

{\bf 11. PROPOSITION : }{\it  Let $s=dq$ be the polar decomposition of a semicircular variable in
some $W^*$-probability space with faithful normal state. Then $(q,d)$ are independent, $q$ is
quarter-circular and the distribution of $d$ is given by $\mu_d(X^{2k})=1$ and
$\mu_d(X^{2k+1})=0$.}

{\bf PROOF: }
Look for instance at the semicircular $(x\mapsto x)\in (L^{\infty}[-1,1],\gamma_{0,1})$.

{\bf 12. THEOREM [2] : }{\it  Let $x=vb$ be the polar decomposition of a circular variable in some
$W^*$-probability space with faithful normal state. Then $v$ is Haar-unitary, $b$ is
quarter-circular and $(v,b)$ is a $*$-free pair.}

{\bf PROOF: }
The theorem is a fairly simple consequence of the proposition 10.

Consider the group $G=\z *(\z\times\z /2\z )$ and denote by $z,t,a$ the images of $1\in \z$,
$(1,\hat{0})\in\z\times(\z /2\z )$ and $(0,\hat{1})\in\z\times(\z /2\z )$ by the canonical
embeddings into $G$. 

Let $u=\lambda_G (z)$, $d=\lambda_G (a)$ and choose a quarter-circular $q\in W^*(\lambda_G(t))$.
Then $(q,d)$ are independent, and the distribution of $d$ is given by $\mu_d(X^{2k})=1$
and $\mu_d(X^{2k+1})=0$ (see (3.4)). By (11), $dq$ is semicircular, so by (10), $c:=udq$ is
circular, and:

- the module of $c$ is $q$, which is a quarter-circular.

- the polar part of $c$ is $ud$, which is obviously a Haar-unitary.

- consider the automorphism $\psi$ of $G$ which is the identity on $\z\times\z /2\z $ and sends
$z\mapsto za$. It extends to a trace-preserving automorphism $\tilde{\psi}$ of $W^*(G)$ which sends
$u\mapsto ud$ and $q\mapsto q$. As $u$ and $q$ are $*$-free, it follows that $ud$ and $q$ are
$*$-free.

{\bf 13. }We give in the end another kind of result which seems to be non-trivial, but which
follows easily by using our formalism:

Let $P\in \c [X]$ be a polynomial. We can write $P={\Sigma}_{j=1}^{j=k}m_jX^{p_j}$ with $p_i\neq
p_j, m_j\neq 0$.

Let $z$ and $n_1...n_k$ be the images of $1\in \z$ and of the $1\in \n$'s by the canonical
embeddings into the free product ${\z}*{\n}^{*k}$; let $\lambda = {\lambda}_{{\z}*{\n}^{*k}}$. By
(7), $\{ zn_1^{p_1},...,zn_k^{p_k},n_1z^{-1},...,n_kz^{-1}\}$ is a prefix, so it is a code, and
$\{\lambda(z)(\lambda (n_j^{p_j})+\lambda (n_j)^*)/2\}_{1\leq j\leq k}$ is a circular system.

Let $a_j$ be complex numbers such that $m_j=a_j^{p_j+1}$. By [1], proposition 2.2., the sum $\Sigma
a_j\lambda (z)(\lambda (n_j)^*+\lambda (n_j)^{p_j})$ is a (non centered) circular variable. By
[1], example 3.4.3., the $R$-transform of $a_j(\lambda (n_j)^*+\lambda (n_j)^{p_j})$ is
$m_jX^{p_j}$. The additivity of the $R$-transform ([1], theorem 3.2.3.), implies that the
$R$-transform of $x=\Sigma a_j(\lambda (n_j)^*+\lambda (n_j)^{p_j})$ is $P$.

Thus, for any polynomial $P$, we can find a random variable $x$ such that:

- the $R$-transform of $x$ is $P$.

- if $u$ is a Haar-unitary $*$-free from $x$, then $ux$ is a (non centered) circular variable.

\vskip3mm{\bf ACKNOWLEDGEMENTS: }I would like to thank G. Skandalis, for directing this work and for
many suggestions and advices; A. Boutet de Monvel, for invinting me in her laboratory during the
acomplishment of this research; E. Germain, for many helpful suggestions; P. Biane, for pointing out
an error in a preliminary version of this paper.
\vskip 1cm
\baselineskip=12pt
{\bf REFERENCES:}
\vskip 5mm
\noindent {\bf [1]} Voiculescu, Dykema, Nica - Free random variables, CRM Monograph Series
 $n^{\circ}1$, AMS
(1993)\\
{\bf [2]}  Voiculescu - Circular and semicircular systems and free product factors,
Progress in Math. 92, Birkh\"auser (1990)\\
{\bf [3]} Lothaire - Combinatorics on Words, Addison-Wesley (1983)
\vskip 1cm

\noindent Universit\'e Paris 7, Aile 45-55, $5^{eme}$ \'etage,\\
2 place Jussieu, 75251 Paris Cedex 05.\\
banica@mathp6.jussieu.fr
\vskip 1cm
\noindent AMS Classification: 46L50

\end{document}